\documentclass[10pt]{amsart}
\usepackage{mathrsfs}
\usepackage{amssymb}

\usepackage[english]{babel}

\newtheorem{thm}{Theorem}[section]

\theoremstyle{definition}
\newtheorem{definition}[thm]{Definition}

\newtheorem{theorem}{Theorem}[section]
\newtheorem{e-proposition}[theorem]{Proposition}

\newtheorem{e-definition}[theorem]{Definition\rm}
\def\og{\leavevmode\raise.3ex\hbox{$\scriptscriptstyle\langle\!\langle$~}}
\def\fg{\leavevmode\raise.3ex\hbox{~$\!\scriptscriptstyle\,\rangle\!\rangle$}}

\def\di{\displaystyle}
\def\f{\frac}

\def\D{\Delta }

\def\e{\varepsilon }

\def\n{\nabla }

\def\om{\omega }

\def\s{\sigma }

\newcommand\EE{\mathbb{E}}

\newcommand\R{\mathbb{R}}
\newcommand\RR{\mathbb{R}}

\newcommand\TT{\mathbb{T}}
\newcommand\id{I\!\!\!I}

\newcommand{\SD}{{\mathcal{D}}}

\newcommand{\cH}{\mathcal{H}}
\newcommand{\cE}{\mathcal{E}}

\newcommand{\dd}{\partial}
\newcommand{\ffrac}[2]{\displaystyle{\frac{#1}{#2}}}
\newcommand{\na}{\nabla}
\newcommand{\oo}{\omega}
\newcommand{\si}{\sigma}

\newcommand{\vep}{\varepsilon}

\newcommand{\ev}{e_t(\vep v)}

\newcommand{\evv}{e_t^{-1}(\vep v)}

\newcommand{\BBox}{\rule{6pt}{6pt}}
\newcommand{\Label}[1]{\label{#1}}
\begin{document}
% place in the next line the header (rubrique) chosen for your article,
% if you know it (you can also have 2, format : Header1/Header2
\centerline{}
%\begin{frontmatter}


% \author{Name\thanksref{label2}}
\selectlanguage{english}
\title{Variational principle for weighted porous media equation}

% The [label1] can be suppressed if there is only one address for all authors

\title[]{Variational principle for weighted porous media equation}
%\date{\today\ \emph{ File: }\jobname.tex}
\author[A. Antoniouk]{Alexandra Antoniouk} \address{Dep. Nonlinear Analysis
\hfill\break\indent Institute of Mathematics NAS Ukraine\hfill\break\indent
 Tereschchenkivska str, 3\hfill\break\indent
  Kyiv, 01 601 UKRAINE}
\email{antoniouk.a@gmail.com}
\author[M. Arnaudon]{Marc Arnaudon} \address{Institut de Mathématiques de Bordeaux \hfill\break\indent CNRS: UMR 5251\hfill\break\indent
  Université Bordeaux 1 \hfill\break\indent
  F33405 TALENCE Cedex, France}
\email{marc.arnaudon@math.u-bordeaux1.fr}

%\selectlanguage{english}
%\author[authorlabel1]{Alexandra Antoniouk},
%%\ead{antoniouk@imath.kiev.ua}
%\author[authorlabel2]{Marc Arnaudon}
%%\ead{marc.arnaudon@math.u-bordeaux1.fr}
%
%\address[authorlabel1]{Dept. Nonlinear Analysis, Institute of Mathematics NAS Ukraine, Tereschchenkivska str. 3,  Kyiv, 01 601 UKRAINE}
%\address[authorlabel2]{Institut de Mathématiques de Bordeaux CNRS: UMR 5251, Universit\'e Bordeaux 1,  F33405, TALENCE Cedex, France}

% If you know the dates of reception, and acceptation you can put them now;
%  idem the name of the person presenting the Note

%\medskip
%\begin{center}
%{\small Received *****; accepted after revision +++++\\
%Presented by £££££}
%\end{center}
\maketitle
\begin{abstract}
\selectlanguage{english}
% Text of abstract in English
In this paper we state the variational principle for the weighted porous media equation. It extends V.I. Arnold's approach to the description of Euler flows as a geodesics on some manifold, i.e. as a critical points of some energy functional.

\vskip 0.5\baselineskip

\end{abstract}
%\end{frontmatter}

% now the Version française abrégée, if it exists
%\selectlanguage{francais}
%\section*{Version fran\c{c}aise abr\'eg\'ee}
% Text of your Version française abrégée here.
% Note you do not need to repeat here equations that you use in the
% main text - for example 'voir (3)' is quite acceptable.

%\selectlanguage{english}
% main text
\section{Introduction}\label{Section1}
In the beginning  18th century Leibniz, Maupertuis, Euler claimed that all physical phenomenons may be obtained from the Least Action Principle and since Lagrange and Hamilton it was well understood for the classical mechanics. However only in
1966 V.I. Arnold in \cite{A} achieved it for the fluid dynamics. To do this he remarked that the group of volume preserving diffeomorphisms $\SD_\mu(M)$ of a  manifold $M$ ($\mu$ being a given volume element on $M$) is the appropriate configuration space for the hydrodynamics of an incompressible fluid. In this framework the solutions to the Euler equation become geodesic curves  with respect to the right invariant metric on $\SD_\mu$, which at $g\in\SD_\mu$ is given by
$
\bigl (X, Y \bigr) = \int_M < X(x), Y(x)>_x d\mu(x),
$
for $X, Y \in T_g \SD_\mu$, $< \cdot ,\cdot >_x$ is a metric on $T_x M$, and $\mu$ is the volume element on $M$ induced by the metric. The relation between geodesics on $\SD_\mu$ and the Euler equation was further studied in \cite{EM} and shortly may be expressed in the following way.
Let $t\mapsto g_t\in \SD_\mu$ be a geodesic with respect to the right invariant metric $(\cdot, \cdot)$, $v_t =\frac{d}{dt}g_t$ be the corresponding velocity, and $u_t = v_t\circ g^{-1}_t$
be a time dependent vector field on M. Then $u_t$ is a solution to the Euler equation for perfect fluid.
%\begin{equation}\Label{Eu}\ 
%\left\{ \displaystyle{
%   \begin{array}{lcl}
%     \ffrac{\dd u_t}{\dd t} +\na_{u_t} u_t   & =    \na p_t,
%\\
%     div\    u_t, \ \ \ u_t\vert_{t=0} &=& 0.
%     \end{array}}
%\right.
%\end{equation}
%Above $p_t: M \mapsto \R$ is the pressure and $\na$ is the covariant derivative. 
In particular the map $t\mapsto g_t$ defined on some time interval $[0,T]$ minimizes the energy functional 
\begin{center}
$\displaystyle 
S(g)=\frac{1}{2}\int_0^T\Big(\int_M \big\|\frac{d g_t}{dt}\big\|^2d\mu(x)\Big)\,dt
$
\end{center}
and the Euler-Lagrange equations for this functional are precisely the Euler equation for perfect fluid. 

Developing this approach in \cite{AC}, \cite{CC}, \cite{Gl}, by means of stochastic methods it was shown that 
an incompressible stochastic flow $g(u)$ with generator $\di \f12\D+u_t$ is critical for some energy functional
%$$
%g\mapsto \f12\EE\left[\int_0^T\int_{\TT} \left\|Dg(t,x)(\om)\right\|^2d\mu(x)\, dt \right], \ \ \mbox{with}\  \ Dg(t,x)(\om) = u(t,g(t,x)(\om))
%$$
if and only if $u$ solves Navier-Stokes equation for viscous incompressible fluid.
% , i.e. there exists a function $P(x)$ such that 
%$$
%\f{\partial u}{\partial t}=-u\cdot \n+\f12\D u+\n P.
%$$
See also~\cite{Constantin-Iyer:08} and \cite{Eyink:10} for other stochastic characterizations of solutions to Navier-Stokes equation.
The purpose of this article is to show that the weighted porous media equation (\cite{DNS}, \cite{DGGW}), which generalizes the standard porous media equation,
\begin{equation}
\Label{PMpq}\ 
\f{\partial u}{\partial t}=\left(-u\cdot \n+\f12\D\right)\left(\|u\|^{q-2}u\right)+\n P.
\end{equation}
may be also obtained in the framework of Least Action Principle for specially chosen energy functional. In the particular case of  $q=2$ this recovers the Navier-Stokes equation.
%\begin{remark}
%\Label{R0}
%For $q=2$ we recover Navier-Stokes equation.
%\end{remark}
%\begin{remark}
%\Label{R1}
%Compare to the porous medium equation
%$$
%\f{\partial u}{\partial t}=\D (\|u\|^{m-1}u).
%$$
%\end{remark}

\section{Operator formulation of variational principle.}\label{Section2}
For simplicity we work on the torus $\TT$ of dimension~$N$. 
%Fix $p\ge 2$. 
From now on, when integrating in the torus, $dx$ will stand for the normalized Lebesgue measure.

\begin{definition} For  some smooth divergence free time dependent vector field $(t,x)\mapsto v_t(x)\in T_x\TT$ we define the flow of $\dot{v}_t$: $e_t(v)\in \SD_\mu(\TT)$  as a solution of the ordinary differential equation
\begin{equation}\Label{pert}\ 
     \ffrac{d e_t(v)}{d t}  =  \dot{v_t}(e_t(v)) ,
\ \ \ 
    e_0(v)
   =
   \id_{\TT}.
\end{equation}
\end{definition}
Let us remark that in some sense $e_t(v)$ is a perturbation of identity map in space $\SD_\mu(\TT)$.
The solvability of this equation easily follows from the compactness of $\TT$ and smoothness of~$v$.

Consider  a time-dependent divergence-free vector field $u$  on $[0,T]\times \TT$. So $u$ takes its values in the tangent bundle of $\TT$ which can at every point be identified with $\RR^N$. Divergence-free means that $\sum_{j=1}^N\partial_ju^j\equiv0$. Define the operator $L(u_t): C^\infty (\TT,\R^N)\to C^\infty (\TT,\R^N)$ by
$ \displaystyle L(u_t)f=\f12\D f+u_t\cdot\n f.
$
%It again takes its values in $\RR^N$.
\begin{definition}
{\it The energy functional} 
 is defined for  $q>1$ as
\begin{equation}\Label{Spq}\ 
\cE_{q}(u,v)=\f1q\int_0^T\int_{\TT}\left\|\Bigl[\big(\partial_t+L(u_t)\big)e_t(v) \Bigr](e_t^{-1}(v)(x))\right\|^q dx\, dt,
\end{equation}
where $e_t^{-1}(v)$ is the inverse map of the diffeomorphism $e_t(v) : \TT\to \TT$.
\end{definition}
\begin{definition}
We say that $u$ is a critical point of $\cE_{q}$  if for all divergence-free time dependent vector field $v$ such that $v_0= 0$ and $v_T= 0$,
$\displaystyle 
\f{d}{d\e}\Big\vert_{\e=0}\cE_{q}(u,\e v)=0.
$
\end{definition}

\begin{theorem}
\Label{P1} 
A  divergence-free time dependent vector field 
$u$  is a critical point of $\cE_{q}$, $q\geq 2$  if and only if there exists a function $P(x)$ such that (\ref{PMpq}) is satisfied.
\end{theorem}
Proof.
\noindent For $
\ev_*\big(u_t\big)(x)=T_{\evv(x)}\ev\left(u_t\left(\evv(x)\right)\right),
$ we compute
\begin{eqnarray*}
\Big[\big(\dd_t+L(u_t)\big)\ev\Big]\big(\evv(x)\big)&=&\vep\dot{v}(t,\evv(x))+\ev_*\big(u_t\big)(x)+\\
&+&
\ffrac{1}{2}\,\big(\Delta \ev\big)\big(\evv(x)\big), 
\end{eqnarray*}
where 
$T_y\ev(\cdot)$ being the tangent map of $\ev$ at point $y$.
Therefore  we have
$$\ffrac{d}{d\vep}\Big\vert_{\vep=0} \Big[\big(\dd_t+L(u_t)\big)
\ev\Big]\big(\evv(x)\big) =
\dot{v}_t(x)+[u_t,v_t](x)+\ffrac12\Delta v_t(x).
$$
Since $u_t=\big(\dd_t+L(u_t)\big)(\id),$
for $\id=e_t(0) : \TT\to\TT$ the identity map,
%$$\SSS_{q,p}\Big(L^p(u), \e v\Big)=\ffrac{1}{p}\int\limits_0^T \int\limits_\TT \Vert \big(\dd_t +L_t^{u,\vep v}\big)(id) \Vert^p dt\, dx
%$$
$\displaystyle \ffrac{d}{d\vep}\Big\vert_{\vep=0} \cE_{q}(u, \e v)$ equals
$$
 \int\limits_0^T \int\limits_\TT\!\! \Vert \big( \dd_t+L(u_t)\big) (id)\Vert^{q-2}
\big<\dot{v}_t+[u_t,v_t]+ \ffrac12\Delta v_t,u_t \big>dx\, dt
=$$
$$=\int\limits_0^T \int\limits_\TT\Vert u_t\Vert^{q-2}\big<\dot{v}_t+[u_t,v_t]+\ffrac12\Delta v_t,u_t \big> \, dx\, dt.
$$
On the other hand
\begin{align*}
&0= \int\limits_\TT \Vert u_T\Vert^{q-2} \big<u_T, v_T\big> dx\\
&= \int\limits_0^T \int\limits_\TT \Big\{ \Vert u_t\Vert^{q-2}\big< u_t,\dot{v}_t\big> + \big<  \Vert u_t\Vert^{q-4}
(q-2)<\dot{u}_t,u_t>u_t + \Vert  u_t \Vert^{q-2}
\dot{u}_t,v_t \big> \Big\}dx\, dt.
\end{align*}
Therefore, writing $u=u_t$ and $v=v_t$,
\begin{align*}
&0= \ffrac{d}{d\vep}\Big\vert_{\vep=0} \cE_{q}(u, \e v) + \\
&+\int\limits_0^T \int\limits_\TT \Big\{\Vert u\Vert^{q-2}\Big(\big<\dot{u}, v\big>- \big< [u,v],u\big> - \ffrac{\big<\D v,u\big>}{2}\Big)+ (q-2) \Vert u\Vert^{q-4} \big< \dot{u},u\big> \big<u,v\big> \Big\} dx\, dt.
\end{align*}
Due to 
$$\displaystyle \int\limits_\TT \Vert u\Vert^{q-2} <\na_v u, u>dx =\ffrac{1}{q}\int\limits_\TT <\na\Vert u\Vert^q, v> dx =-\ffrac{1}{q}\int\limits_\TT \Vert u\Vert^q \, \mbox{div}  \, v\, dx = 0
$$
for $\mbox{div}\, v =0$, we have, using $[u,v]=\n_uv-\n_vu$,
\begin{eqnarray*}
-\ffrac{d}{d\vep}\Big\vert_{\vep =0} \cE_{q}\Big( u,(\vep v)\Big) &=& 
\int\limits_0^T \int\limits_\TT \Big\{-\Vert u\Vert^{q-2} \big< \na_u 
v,u\big>-\ffrac12\big<v, \Delta \Big( \Vert u\Vert^{q-2}u\Big)\big>+\\
&+&(q-2)\Vert u\Vert^{q-4}\big< \dot{u},u\big>\big<u,v\big>+ \Vert u
\Vert^{q-2}\big<\dot{u}, v\big> \Big\} \, dx\, dt\\
&=& \int\limits_0^T \int\limits_\TT \Big< \na_u\Big( \Vert u\Vert^{q-2}u\Big)-
\ffrac12\Delta \Big( \Vert u\Vert^{q-2} u \Big)+\\
&+&(q-2)\Vert u\Vert^{q-4}\big<\dot{u},u\big>u+ \Vert u\Vert^{q-2}\dot{u}, v\Big> \, dx\, dt=\\
&=& \int\limits_0^T\int\limits_\TT \big< \big(\dd_t +u\cdot \na -\ffrac12 \Delta\big)\Vert u\Vert^{q-2}u,v\big> \, dx\, dt
\end{eqnarray*}
(notice that in the second equality we used the fact that $\int_\TT u\left(\langle v,\|u\|^{q-2}\rangle\right) dx=\int_\TT \mbox{div}\   u\langle v,\|u\|^{q-2}\rangle dx=0$).
This equality is true for all time dependent divergence free vector field $v$, so it gives the equivalence between $u$ critical point of $\cE_{q}$ and solution to equation~(\ref{PMpq}).
\BBox

\section{Stochastic variational principle for incompressible diffusion flows}\label{Section3}
%\setcounter{equation}0

%Now we consider any compact oriented Riemannian manifold $M$ without boundary of dimension $N\ge 2$.  We shall use notation $dx$ for integration with respect to the normalized volume measure on $M$.

We define a  {\it diffusion flow} $g_{{t}}(x)$ on $\TT$, $x\in \TT$,  $t\in[0,T]$, $T>0$ as a stochastic process, which satisfies
 the It\^o {stochastic} equation:
\begin{equation}\Label{1}\ 
dg_t(x)=\s(g_t(x))\,dW_t+u_t(g_t(x))\,dt,\quad g_0(x)=x
\end{equation}
where $u_t$ is a time dependent vector field on $\TT$, $\si\in\Gamma(Hom(\cH, T\TT))$ is a $C^2$-map satisfying for all $x\in \TT$ $(\si\si^\ast)(x)=\id_{T_x\TT}$, $W_t$ is a cylindric Brownian motion in Hilbert space $\cH$. 

Let us remark that a diffusion flow is a diffusion process ${\{}g_t(u)(x){\}}_{t\ge 0}$ with generator $L(u_{{t}}) =\f12\D+u_{{t}}$.
We define an {\it incompressible 
 diffusion flow}
$g_t(u)(x)(\om)$ as a diffusion flow such that a.s. $\oo$ for all $t\ge 0$, the map $x\mapsto g_t(u)(x)(\om)$ is a volume preserving diffeomorphism of~$\TT$.
Examples of incompressible diffusion flows can be found in~\cite{CC}. Notice that a necessary condition is $\mbox{div}\,  u_t = 0$.

For the diffusion flow $g_t$ (\ref{1}) 
%of the form
%$$
%g_t=g_0+\int_0^t {\si}_s(g_s) dW_s+\int_0^t  u_s(g_s,\oo) ds,
%$$ 
we define the {\it drift} as the time derivative of the finite variation part by 
$
D g_t (\oo):=  u_t(g_t,\oo),
$ and 
 the {\it energy functional} {by}
\begin{equation}\Label{Eq}\ 
\cE_q(g):=\f1q\EE\Big[\int_0^T\int_{\TT} \Big\|Dg_t(x)(\om)\Big\|^q\, dx\, dt\Big], \ \ \ q>1. 
\end{equation}

We make a perturbation  by letting $g_t^v(u)=e_t(v)\circ g_t(u)$, where {$v$} is a smooth divergence free time dependent vector field and $e_t(v)$ is defined in~(\ref{pert}). 

\begin{definition}{\Label{D3.2}}\ 
We say that $g_t(u)$ is a critical point for the energy functional $\cE_q$ if for all smooth time dependent divergence free vector field  $v$ on $T\TT$ such that $v_0=v_T=0$, 
$$
\f{d }{d\e}\Big\vert_{\e=0}\cE_q(g^{\e v}(u))=0.
$$
\end{definition}

\begin{thm}
\label{C1}
Let $q\ge  2$.
An incompressible diffusion flow $g_{{t}}(u)$ with generator $\di L(u_{{t}})$ is {a} critical  {point} for the energy functional $\cE_q$
if and only if there exists a function $P(x)$ such that $u_t$ satisfies equation (\ref{PMpq}).
%$$
%\f{\partial u}{\partial t}=\left(-u\cdot \n+\f12\D\right)\left(\|u\|^{q-2}u\right)+\n P.
%$$
\end{thm}
{\bf Proof} of this theorem is a consequence of Theorem~\ref{P1}  and the  It\^o's formula. 
%Applying It\^o's formula we get 
%$$
%Dg_t^{\e v}(u)(x)=\Big[\left(\partial_t+L(u_t)\right)(e_t(\e v))\Big]\left(e_t(\e v)^{-1}(g_t^{\e v}(u)(x))\right)
%$$
%so we have 
%$$
%\cE_q(g^{\e v}(u))=\f1q\EE\,\int_0^T\int_{\TT} \left\|\Big[\left(\partial_t+L(u_{{t}})\right)(e_t(\e v))\Big]\left(e_t(\e v)^{-1}(g_t^{\e v}(u)(x))\right)\right\|^q\, dx \, dt.
%$$
%But  {for} a.s.  {$\oo$ and} all $t\in[0,T]$ the map $x\mapsto g_t^{\e v}(u)(x)(\om)$ is measure preserving. Thus making the change of variable $g_t^{\e v}(u)(x)(\om)=y$ we can remove the expectation in the formula {above} and we get $\cE_q(g^{\e v}(u))=\cE_{q}(u,\e v)$. So the rest of the proof  immediately follows from Theorem~\ref{P1}.
\BBox

%Similarly we have the following result.
%\begin{thm}
%\label{C2}
%The  incompressible  stochastic flow $g_t(u)$ with generator $\di L(u_t)$ is a critical point for the energy functional
%$$
%g\mapsto  {\cE^1_q(g)=}\f1q\EE\left[\int_0^Tdt\int_{\TT}dx \left\|\n Dg_t(u)(x)(\om)\right\|^q\right]
%$$
%if and only if there exists a function $P(x)$ such that 
%$$
%\f{\partial \D_q u}{\partial t}=\left(-u\cdot \n-(\n u)^\ast+\f12\D\right)\left(\D_qu\right)+\n P.
%$$
%\end{thm}
%The proof follows the same lines as the one of Theorem~\ref{C1} and is a consequence of Theorem~\ref{P3}
% etc, etc

% \appendix

% \section{}
% \label{}

% The Acknowledgements are an un-numbered section
\section*{Acknowledgements}
% Acknowledgements text here
The research of A.Antoniouk was supported by the grant no. 01-01-12
of National Academy of Sciences of Ukraine (under the joint
Ukrainian-Russian project of NAS of Ukraine and Russian
Foundation of Basic Research)."

This research was also supported by the French ANR grant ANR-09-BLAN-0364-01 ProbaGeo.

\end{document}